# An Integral over (0, π) for the Distribution Function of a Sum of Independent Gamma Random Variables and for Quadratic Forms of Gaussian Variables

Thomas Royen

TH Bingen, University of Applied Sciences

Berlinstrasse 109, D55411 Bingen, Germany

e-mail: thomas.royen@t-online.de

**Abstract.**

An integral over the interval $(0, \pi)$ is given for the cumulative distribution function (cdf) of a sum of independent gamma random variables with different scale and shape parameters. The cdf of a positive definite quadratic form is obtained as a special case with identical shape parameters $\alpha = \frac{1}{2}$.

## 1. Introduction

There is a lot of literature on the distribution of quadratic forms, in particular with Gaussian random vectors, see e.g. the extensive chapter 4 in [3]. For all these series representations of the cdf a large number of coefficients has to be computed. It is the purpose of the underlying paper to avoid these computations and to use only the eigenvalues of a positive definite matrix.

Let $X$ be an $N_k(0, \Sigma)$-random vector with a pos. def. covariance matrix $\Sigma$. For the quadratic form $Q = X'CX$ with a pos. def. matrix $C$ we have the equation

$$P\{Q \le x\} = P\left\{\sum_{j=1}^{k} \lambda_j Z_j^2 \le x\right\} \tag{1.1}$$

with the eigenvalues $\lambda_j$ of $\Sigma^{1/2} C \Sigma^{1/2}$ and an $N_k(0, I_k)$- random vector $Z$, where $I_k$ denotes the $k \times k$ - identity matrix. The determinant of a matrix $M$ will be denoted by $|M|$, and $\|M\|$ is the spectral norm of a symmetrical matrix.

The random variable $\frac{1}{2} Z_j^2$ has a gamma-pdf $g_\alpha(x) = \exp(-x) x^{\alpha-1} / \Gamma(\alpha)$ with $\alpha = \frac{1}{2}$. Therefore, the distribution of $\frac{1}{2} Q$ is a special case of the distributions of sums of independent gamma random variables with pdfs $\lambda_j^{-1} g_{\alpha_j}(\lambda_j^{-1} x)$. The pdf $f(x; \alpha_1, ..., \alpha_k, \lambda_1, ..., \lambda_k)$ of such a sum has the Laplace transform (Lt)

$$\prod_{j=1}^{k} (1 + \lambda_j t)^{-\alpha_j}. \tag{1.2}$$

For series expansions of the corresponding cdfs see [4] and for multivariate generalizations of these series see [5].

---

## 2. The Integral Representations

For the intended integral representations we use a modified form of the Lt from (1.2). With

$$v = \tfrac{1}{2}\left(\lambda_{max}^{-1} + \lambda_{min}^{-1}\right), \ z = (1+v^{-1}t)^{-1}, \text{ and } \alpha = \sum_{j=1}^{k} \alpha_j \tag{2.1}$$

we obtain

$$\prod_{j=1}^{k}(1+\lambda_j t)^{-\alpha_j} = \prod_{j=1}^{k}(1+\lambda_j v t v^{-1})^{-\alpha_j} = \prod_{j=1}^{k}(1 + v^{-1}t + (\lambda_j v - 1)tv^{-1})^{-\alpha_j} =$$

$$z^{\alpha}\prod_{j=1}^{k}(1+(\lambda_j v - 1)(1-z))^{-\alpha_j} = \left(v^{-\alpha}\prod_{j=1}^{k}\lambda_j^{-\alpha_j}\right)z^{\alpha}\prod_{j=1}^{k}\left(1-(1-\tfrac{1}{v\lambda_j})z\right)^{-\alpha_j}, \tag{2.2}$$

and $\max |c_j| = \max\left|1 - \dfrac{1}{v\lambda_j}\right| = \dfrac{\lambda_{max} - \lambda_{min}}{\lambda_{max} + \lambda_{min}} < 1.$

To avoid series expansions with products of binomial series the orthogonality of the functions $e^{in\varphi}$, $n \in \mathbb{Z}$, will be used for a product of a power series of $y = re^{i\varphi}$ and a power series of $r^{-1}e^{-i\varphi}$. For this goal the following function is defined for $\alpha > 0$, $x \geq 0$, $y \in \mathbb{C}$:

$$G_{\alpha}(x,y) := \sum_{n=0}^{\infty} G_{\alpha+n}(x) y^n \tag{2.3}$$

where $G_{\alpha}(x)$ is the cdf with the gamma-pdf $g_{\alpha}(x) = e^{-x}x^{\alpha-1}/\Gamma(\alpha)$. Here we need only values $|y| < 1$.

Now the equation

$$G_{\alpha}(x,y) = (1-y)^{-1}\left(G_{\alpha}(x) - y^{1-\alpha}e^{(y-1)x}G_{\alpha}(xy)\right). \tag{2.4}$$

is verified. Integration by parts shows $G_{\alpha}(y) = \sum_{n=1}^{\infty} g_{\alpha+n}(y)$, and it is

$$(1-y)G_{\alpha}(x,y) = G_{\alpha}(x) + \sum_{n=1}^{\infty}(G_{\alpha+n}(x) - G_{\alpha+n-1}(x))y^n = G_{\alpha}(x) - \sum_{n=1}^{\infty} g_{\alpha+n}(x)y^n =$$

$$G_{\alpha}(x) - y^{1-\alpha}e^{(y-1)x}\sum_{n=1}^{\infty} g_{\alpha+n}(xy) = G_{\alpha}(x) - y^{1-\alpha}e^{(y-1)x}G_{\alpha}(xy).$$

For $y = 1$ we obtain $G(x,1) := \lim_{y\to 1} G(x,y) = xg_{\alpha}(x) + (1+x-\alpha)G_{\alpha}(x)$. Furthermore, a short calculation shows

$$G_{\alpha}(x,y) = (1-y)^{-1}\left(G_{\alpha-1}(x) - y^{1-\alpha}e^{(y-1)x}G_{\alpha-1}(xy)\right), \ \alpha \geq 1, \ G_0 := 1. \tag{2.5}$$

In particular, it is $G_{1+n}(y) = 1 - e^{-y}\sum_{k=0}^{n}\dfrac{y^k}{k!}$ and $G_{1/2+n}(y) = erf(y^{1/2}) - e^{-y}\sum_{k=1}^{n}\dfrac{y^{k-1/2}}{\Gamma(k+1/2)}$.

With binomial series we obtain for the Lt from (2,2)

$$\prod_{j=1}^{k}(1+\lambda_j t)^{-\alpha_j} = \left(v^{-\alpha}\prod_{j=1}^{k}\lambda_j^{-\alpha_j}\right)\sum_{n=0}^{\infty}\left(\sum_{(n)}\prod_{j=1}^{k}\binom{\alpha_j + n_j - 1}{n_j}c_j^{n_j}\right)z^{\alpha+n}, \tag{2.6}$$

where $\Sigma_{(n)}$ stands for $\sum_{n_1+\ldots+n_k=n}$ and consequently for the corresponding cdf:

$$F(x;\alpha_1,\ldots,\alpha_k,\lambda_1,\ldots,\lambda_k) = \left(v^{-\alpha}\prod_{j=1}^{k}\lambda_j^{-\alpha_j}\right)\sum_{n=0}^{\infty}\left(\sum_{(n)}\prod_{j=1}^{k}\binom{\alpha_j + n_j - 1}{n_j}c_j^{n_j}\right)G_{\alpha+n}(vx). \tag{2.7}$$

Then with $G_{\alpha}(x,y)$ from (2.4) or (2.5) it follows directly:

**Theorem.** Let $X_1,...,X_k$ be independent gamma random variables with pdfs $\lambda_j^{-1} g_{\alpha_j}(\lambda_j^{-1} x)$. Then the cdf of $\sum_{j=1}^{k} X_j$ is given by

$$F(x;\alpha_1,...,\alpha_k,\lambda_1,...,\lambda_k) = \left(v^{-\alpha} \prod_{j=1}^{k} \lambda_j^{-\alpha_j}\right) \frac{1}{\pi} \int_0^\pi Re\left(\prod_{j=1}^{k}(1-c_j r^{-1} e^{-i\varphi})^{-\alpha_j}\right) G_\alpha(vx, re^{i\varphi}) d\varphi \quad (2.8)$$

with $v$ from (2.1), $c_j = 1 - \frac{1}{v\lambda_j}$ from (2.2) and any $r \in \left(\frac{\lambda_{\max} - \lambda_{\min}}{\lambda_{\max} + \lambda_{\min}}, 1\right)$.

In particular, the cdf of the quadratic form $\frac{1}{2}Q$ from (1.1) with an $N_k(0,\Sigma)$ - random vector $X$ is given by

$$F(x; \tfrac{1}{2}, \lambda_1,...,\lambda_k) = \left(v^{-k/2} |C\Sigma|^{-1/2}\right) \frac{1}{\pi} \int_0^\pi Re\left(\prod_{j=1}^{k}(1-c_j r^{-1} e^{-i\varphi})^{-1/2}\right) G_{k/2}(vx, re^{i\varphi}) d\varphi. \quad (2.9)$$

For multivariate generalizations of formula (2.8) see [5], but the numerical evaluation of the corresponding multivariate integrals is frequently difficult because of a rather high total variation of the integrands. Here only such a formula is given for the cdf of a single p-variate gamma distribution with the Lt $|I_p + \Sigma T|^{-\alpha}$, where $\Sigma$ is pos. def. with the eigenvalues $\lambda_k$ and $T := Diag(t_1,...,t_p)$, $t_k \geq 0$. This distribution exists at least for $2\alpha \in \mathbb{N} \cup ([\frac{p-1}{2}], \infty)$, see [6], and for infinitely divisible multivariate gamma distributions for all $\alpha > 0$, see [1]. The corresponding cdf is

$$F(x_1,...,x_p; \alpha, \Sigma) = \left(v^{-\alpha p} |\Sigma|^{-\alpha}\right) |I_p - CY^{-1}|^{-\alpha} \frac{1}{(2\pi)^p} \int_{(-\pi,\pi)^p} \prod_{k=1}^{p} G_\alpha(vx_k, y_k) d\varphi_k \quad (2.10)$$

with $Y = Diag(y_1,...,y_p)$, $y_k = re^{i\varphi_k}$, $C = I_p - (v\Sigma)^{-1}$, $\|C\| = \frac{\lambda_{\max} - \lambda_{\min}}{\lambda_{\max} + \lambda_{\min}} < 1$ with the eigenvalues $\lambda_k$ of the pos. def. matrix $\Sigma$ and any $r \in \left(\frac{\lambda_{\max} - \lambda_{\min}}{\lambda_{\max} + \lambda_{\min}}, 1\right)$. However, there are numerical more favourable integral representations for these multivariate gamma cdfs with $p \leq 4$, see [2]. The appendix of this paper contains presumably the most comprehensive collection of formulas for multivariate gamma distributions of this type.